\documentclass[12pt]{article}
\usepackage{latexsym}
\usepackage{amsmath}
\usepackage{amssymb}
\usepackage{amsfonts}
\usepackage{tikz}

\numberwithin{equation}{section}
\newtheorem{theorem}{Theorem}[section]
\newtheorem{lemma}[theorem]{Lemma}

\newtheorem{claim}[theorem]{Claim}

\newtheorem{remark}[theorem]{Remark}
\newcommand{\eproof}{{\mbox{\ }~\hfill
\mbox{\large $\Box$} \par \vskip 10pt}}

\newcommand{\supp}{\mbox{\rm supp\,}}

\newcommand{\eps}{\varepsilon}
\newcommand{\R}{{\mathbf R}}

\newcommand{\pf}{\noindent{\bf Proof}.\;}

\renewcommand{\d}{\partial}
\newcommand{\vphi}{\varphi}
\newcommand{\bd}{\bar\partial}

\title{On Landis' conjecture in the plane}
\author{Carlos Kenig\thanks{Department of Mathematics, University of Chicago, Chicago, IL 60637, USA. Email: cek@math.uchicago.edu. Supported in part by NSF Grant DMS-1265249.}\qquad Luis Silvestre\thanks{Department of Mathematics, University of Chicago, Chicago, IL 60637, USA. Email: luis@math.uchicago.edu. Supported in part by NSF grants DMS-1254332 and DMS-1065979.}\qquad Jenn-Nan Wang\thanks{Institute of Applied Mathematical Sciences, NCTS (Taipei), National Taiwan University,
Taipei 106, Taiwan. Email: jnwang@math.ntu.edu.tw. Supported in part by the Ministry of Science and Technology. }}

\date{}

\begin{document}
\maketitle

\begin{abstract}
In this paper we prove a quantitative form of Landis' conjecture in the plane. Precisely, let $W(z)$ be a measurable real vector-valued function and $V(z)\ge 0$ be a real measurable scalar function, satisfying $\|W\|_{L^{\infty}(\R^2)}\le 1$ and $\|V\|_{L^{\infty}(\R^2)}\le 1$. Let $u$ be a real solution of $\Delta u-\nabla(Wu)-Vu=0$ in $\R^2$. Assume that $u(0)=1$ and $|u(z)|\le\exp(C_0|z|)$. Then $u$ satisfies
$\underset{{|z_0|=R}}{\inf}\,\underset{|z-z_0|<1}{\sup}|u(z)|\ge \exp(-CR\log R)$, where $C$ depends on $C_0$. In addition to the case of the whole plane, we also establish a quantitative form of Landis' conjecture  defined in an exterior domain.
\end{abstract}

\section{Introduction}\label{sec1}
\setcounter{equation}{0}

In the late 60's (see \cite{kl88}), E.M. Landis  conjectured that if $\Delta u+Vu=0$ in $\R^n$ with $\|V\|_{L^{\infty}(\R^n)}\le 1$ and $\|u\|_{L^{\infty}(\R^n)}\le C_0$ satisfying $|u(x)|\le C\exp(-C|x|^{1+})$, then $u\equiv 0$. Landis' conjecture was disproved by Meshkov \cite{m92} who constructed such $V$ and nontrivial $u$ satisfying $|u(x)|\le C\exp(-C|x|^{\frac 43})$. He also showed that if $|u(x)|\le C\exp(-C|x|^{\frac 43+})$, then $u\equiv 0$. A quantitative form of Meshkov's result was derived by J. Bourgain and the first author \cite{bk05} in their resolution of Anderson localization for the Bernoulli model \cite{an58} in higher dimensions. It should be noted that both $V$ and $u$ constructed by Meshkov are \emph{complex-valued} functions. It remains an open question whether Landis' conjecture is true for real-valued $V$ and $u$. In this paper, we confirmed Landis' conjecture for any real solution $u$ of
\begin{equation}\label{e1}
\Delta u-\nabla(Wu)-Vu=0\quad\text{in}\quad \R^2,
\end{equation}
where $W(x,y)=(W_1(x,y),W_2(x,y))$ and $V(x,y)$ are measurable, real-valued, and $V(x,y)\ge 0$ a.e.  In view of the scaling argument in \cite{bk05}, Landis' conjecture is an easy consequence of the estimate for the maximal vanishing order of $u$ satisfying \eqref{e1} in a bounded domain (see also \cite{ke07}). To be precise, we prove that
\begin{theorem}\label{t1}
Assume that  $W_1(x,y), W_2(x,y)$ and $V(x,y)$ are measurable, real-valued, and $V(x,y)\ge 0\ {a.e.}$ in $B_2$, moreover, there exist $K\ge 1$, $M\ge 1$ such that
\[
\|W\|_{L^{\infty}(B_2)}\le K,\quad\|V\|_{L^{\infty}(B_2)}\le M.
\]
Let $u$ be a real solution to
\begin{equation}\label{e2}
\Delta u-\nabla(Wu)-Vu=0\quad\text{in}\quad B_2.
\end{equation}
Assume that $\|u\|_{L^{\infty}(B_2)}\le\exp(C_0(\sqrt{M}+K))$ and 
\[
\|u\|_{L^{\infty}(B_{1})}\ge 1.
\]
Then 
\begin{equation}\label{r1}
\|u\|_{L^{\infty}(B_r)}\ge r^{C(\sqrt{M}+K)}
\end{equation}
for all sufficiently small $r$, where $C$ depends on $C_0$.
\end{theorem}

Hereafter, we denote $B_r(a)$ an open disc of radius $r$ centered at $a$. In the case when $a=0$, we simply denote $B_r(0)=B_r$. Having proved Theorem~\ref{t1}, Landis' conjecture in quantitative form follows by a scaling argument \cite{bk05}.
\begin{theorem}\label{t2}
Assume that  $W_1(x,y), W_2(x,y)$ and $V(x,y)$ are measurable, real-valued, and $V(x,y)\ge 0\ {a.e.}$ in $\R^2$, furthermore, 
\begin{equation}\label{c1}
\|W\|_{L^{\infty}(\R^2)}\le 1,\quad\|V\|_{L^{\infty}(\R^2)}\le 1.
\end{equation}
Let $u$ be a real solution to \eqref{e1}. Assume that $|u(z)|\le\exp(C_0|z|)$ and $u(0)=1$, where $z=(x,y)$. Let $z_0=(x_0,y_0)$. Then we have that
\begin{equation}\label{r2}
\inf_{|z_0|=R}\sup_{|z-z_0|<1}|u(z)|\ge \exp(-CR\log R)\quad\text{for}\quad R\gg 1
\end{equation}
where $C$ depends on $C_0$.
\end{theorem}

Using the same techniques, we can also confirm Landis' conjecture for 
\begin{equation}\label{1-101}
\Delta u+W\cdot\nabla u-Vu=0\quad\text{in}\quad\R^2
\end{equation}
with $V\ge 0$. As before, we first prove an estimate of the maximal vanishing order of $u$ in $B_2$.
\begin{theorem}\label{t3}
Assume that $W$ and $V$ satisfy the same assumptions as in Theorem~\ref{t1}. Let $u$ be a real solution to
\[
\Delta u+W\cdot\nabla u-Vu=0\quad\text{in}\quad B_2.
\]
Assume that $\|u\|_{L^{\infty}(B_2)}\le\exp(C_0(\sqrt{M}+K))$ and 
\[
\|u\|_{L^{\infty}(B_{1})}\ge 1.
\]
Then 
\begin{equation}\label{r18}
\|u\|_{L^{\infty}(B_r)}\ge r^{C(\sqrt{M}+K)}
\end{equation}
for all sufficiently small $r$, where $C$ depends on $C_0$.
\end{theorem}
By Bourgain-Kenig's scaling argument, we then show that
\begin{theorem}\label{t4}
Assume that $W, V$ satisfies the assumptions described in Theorem~\ref{t2}. Let $u$ be a real solution to \eqref{1-101}. Assume that $|u(z)|\le\exp(C_0|z|)$ and $u(0)=1$. Then we have that
\begin{equation}\label{r21}
\inf_{|z_0|=R}\sup_{|z-z_0|<1}|u(z)|\ge \exp(-CR\log R)\quad\text{for}\quad R\gg 1
\end{equation}
where $C$ depends on $C_0$.
\end{theorem}

In this paper we also establish a quantitative form of Landis' conjecture in an exterior domain. Precisely, we show that
\begin{theorem}\label{t5}
Let $u$ be a real solution of
\[
\Delta u-V(x,y)u=0\quad\text{in}\quad B_1^c,
\]
where $V(x,y)\ge 0\ {a.e.}$ is measurable, real-valued, and satisfies 
\[
\|V\|_{L^{\infty}(B_1^c)}\le 1.
\]
Assume that $\|u\|_{L^{\infty}(B_1^c)}\le 1$ and there exists $C_0>0$ such that
\begin{equation}\label{apriorib}
\inf_{|z_0|=\frac 52}\int_{B_1(z_0)}|u|^2\ge C_0.
\end{equation}
Then we have that
\begin{equation*}
\inf_{|z_0|=R}\sup_{|z-z_0|<1}|u(z)|\ge C \exp(-C'R(\log R)^2)\quad\text{for}\quad R\gg 1,
\end{equation*}
where $C$ depends on $C_0$ and $C'$ is an absolute constant.
\end{theorem}

Before outlining the ideas of our proof, we remark on some related works. The exponent $\sqrt{M}+K$ in \eqref{r1} is known to be optimal. For the case where $u$ is a $\lambda$-eigenfunction of the Laplace-Beltrami operator in a smooth compact Riemannian manifold without boundary, the maximal vanishing order of $u$ is less than ${C\sqrt{\lambda}}$ proved by Donnelly and Fefferman in \cite{df88}. Donnelly and Fefferman's proof was based on the Carleman method. Using the method of the frequency function developed by Garofalo and Lin \cite{gl86, gl87}, Kukavica \cite{ku98} proved the maximal vanishing order of $u$ solving $L u+Vu=0$ in $\Omega\subset\R^n, n\ge 2$ is less than $C(1+\|V_-\|^{1/2}_{L^{\infty}(\Omega)}+(\text{osc}_{\Omega}V)^2)$, where $L$ is a general uniform second order elliptic operator, $V_-=\max\{-V,0\}$, and $\text{osc}_{\Omega}V=\sup_{\Omega}V-\inf_{\Omega}V$. However, for the equation \eqref{e2} with $W\not\equiv 0$, Kukavica's method can not produce an order which is algebraic in $\|W\|_{L^{\infty}}$ (see \cite[Remark~5.4]{ku98}).  For \eqref{e2} without the presence of $W$, the method in \cite{bk05} proves that the maximal vanishing order of $u$ is less than $CM^{2/3}$ (see \cite{ke07}). In \cite{da12}, Davey used the Carleman method to study the quantitative uniqueness estimate of $u$ to 
\[
\Delta u+W\cdot\nabla u+Vu=0\quad\text{in}\quad\R^n,
\]
where $W$ and $V$ also satisfy some decaying properties. Her results lead to an order $R^2$ in \eqref{r2} under the assumption \eqref{c1} (see also \cite{lw13}). Moreover, Meshkov-type examples are constructed in \cite{da12} showing that in the presence of $W$ (complex-valued) with $V\equiv 0$, the exponent $2$ is optimal. 

The Carleman method has been a powerful technique in studying questions related to Landis' conjecture and is able to produce optimal bounds in the complex-valued case. Since the Carleman estimate does not seem to distinguish real- or complex-valued functions, a direct use of such estimates to resolve Landis' conjecture seems likely to fail. In this paper we will take a different approach. The main idea lies in the nice relation between second order elliptic equations in the plane and the Beltrami system. The assumption of $V\ge 0$ allows us to construct a global positive multiplier which enables us to convert \eqref{e2} into an elliptic equation in divergence form (see \eqref{3}). By the Beltrami system, we can derive three-ball inequalities (with precise exponent, constant, and radii) for solutions $u$ of \eqref{e2} without using the Carleman estimate or the frequency function.

For the problem in an exterior domain, the local result like Theorem~\ref{t1} remains true. Unfortunately, the usual scaling argument does not lead us to Landis' conjecture.  One of the obstacles is the lack of simply-connectedness in $B_2$ in the rescaled problem. To overcome the difficulty, we introduce an appropriate cutoff function and construct an "approximate" stream function. We reduce the original equation to an inhomogeneous d-bar equation and then apply a Carleman estimate for $\bd$.   

The paper is organized as follows. In Section~\ref{sec2}, we consider a simple $\bd$ equation and study \eqref{e2} with $W\equiv 0$ to explain the main ideas of the proof. Based on the ideas in Section~\ref{sec2}, we prove Theorem~\ref{t1} and \ref{t2} in Section~\ref{sec3}. We prove Theorem~\ref{t3}, \ref{t4} in Section~\ref{sec4}. We prove Theorem~\ref{t5} in Section~\ref{sec5}.  Finally, we discuss some related questions in Section~\ref{conl}. Throughout the paper, $C$ denotes an absolute positive constant whose dependence will be specified whenever necessary. The value of $C$ may vary from line to line. Also, in the paper, we use $z\in\R^2$ or $(x,y)\in\R^2$.

\section{Main ideas and proofs of Theorem~\ref{t1}, \ref{t2} for $W\equiv 0$}\label{sec2}

To motivate the main ideas of our method. We begin with a simple $\bd$ equation.
\begin{equation}\label{equ2.1}
\bd u=V(z)u\quad\text{in}\quad B_2,
\end{equation}
with $\|V\|_{L^{\infty}(B_2)}\le M$. Assume that $\|u\|_{L^{\infty}(B_1)}\ge 1$ and $\|u\|_{L^{\infty}(B_2)}\le e^{C_0M}$. It is clear that any solution of \eqref{equ2.1} can written as
\[
u=\exp(w)f\quad\text{or}\quad f=\exp(-w) u,
\]
where
\begin{equation}\label{ct}
w(z)=-\frac{1}{\pi}\int_{B_2}\frac{V(\zeta)}{\zeta-z}d\zeta.
\end{equation}
Recall that 
\[
|w(z)|\le C\|V\|_{L^{\infty}(B_2)}\le CM\quad\text{for}\quad z\in B_2.
\]
We can see that
\[
\|u\|_{L^{\infty}(B_1)}\ge e^{-CM}\quad\text{and}\quad\|u\|_{L^{\infty}(B_2)}\le e^{CM},
\]
where the second constant $C$ depends on $C_0$.

Since $f$ is holomorphic in $B_2$, by Hadamard's three-circle theorem, we have
\[
\|f\|_{L^{\infty}(B_{r_1})} \le \|f\|_{L^{\infty}(B_{r})}^{\theta}\|f\|_{L^{\infty}(B_{r_2})}^{1-\theta},
\]
where $r<r_1<r_2$ and
\[
\theta=\frac{\log(\frac{r_2}{r_1})}{\log(\frac{r_2}{r})}.
\]
Taking $r_1=1$ and $r_2=\frac 32$ yields
\[
e^{-CM}\le \|u\|_{L^{\infty}(B_{r})}^{\theta}
\]
and hence
\begin{equation}\label{2200}
\|u\|_{L^{\infty}(B_{r})}\ge r^{CM}
\end{equation}
and $C$ depends on $C_0$.

Having derived the estimate of vanishing order, we can prove the following quantitative uniqueness estimate.
\begin{theorem}\label{t0}
Let $u$ be any solution of 
\[
\bd u=Vu\quad\text{in}\quad\R^2.
\]
Assume that $\|V\|_{L^{\infty}(\R^2)}\le 1$, $|u(z)|\le e^{C_0|z|}$, and $u(0)=1$. Then 
\begin{equation}\label{t0e}
\inf_{|z_0|=R}\sup_{|z-z_0|<1}|u(z)|\ge \exp(-CR\log R)\quad\text{for}\quad R\gg 1
\end{equation}
where $C$ depends on $C_0$.
\end{theorem}
\pf This theorem is an easy consequence of \eqref{2200} by the scaling argument used in \cite{bk05}. We include the proof for the sake of completeness. Denote $|z_0|=R$. Let $u_R(z)=u(R(z+z_0/R))$, then $u_R$ satisfies \eqref{equ2.1} with a rescaled potential $V_R(z)=RV(R(z+z_0/R))$, where
\[
\|q_R\|_{L^{\infty}(B_2)}\le R.
\]
It is easy to see that $|u_R(z)|\le e^{C_0R}$ and
\[
u_R(-\frac{z_0}{R})=u(0)=1
\]
with $|\frac{z_0}{R}|=1$. Hence, we have
\[
\|u_R\|_{L^{\infty}(B_{1})}\ge 1.
\]
The quantitative estimate \eqref{t0e} follows easily from \eqref{2200} by taking $r=R^{-1}$.\eproof

We now would like to prove Theorems~\ref{t1} and \ref{t2} without the presence of $W$ based on the ideas described above. We first consider the local problem. Let $u$ be a real solution to
\begin{equation}\label{1}
\Delta u-V(z)u=0\quad\text{in}\quad B_2\subset\R^2,
\end{equation}
where $V(z)\ge 0\ \text{a.e.}$ Likewise, we assume that $u$ satisfies
\begin{equation}\label{1-2}
\|u\|_{L^{\infty}(B_2)}\le \exp(C_0\sqrt{M})
\end{equation}
for some $C_0>0$ and 
\begin{equation}\label{1-1} 
\|V\|_{L^{\infty}(B_2)}\le M. 
\end{equation}
We let $M\ge 1$ for simplicity. We will first construct a positive multiplier. Let us consider $\phi_1(z)=\phi_1(x,y)=\exp(\lambda x)$, then $\Delta\phi_1=\lambda^2\phi$. Hence, if we choose $\lambda=\sqrt{M}$, then we have
\[
\Delta\phi_1-V\phi_1=(\lambda^2-V)\phi_1\ge 0,
\]
that is, $\phi_1$ is a subsolution. On the other hand, define $\phi_2=\exp(2\sqrt{M})$, then $\Delta\phi_2-V\phi_2\le 0$, i.e., $\phi_2$ is a supersolution. Note that $\phi_2\ge\phi_1$. Thus, there exists a positive solution $\phi$ satisfying \eqref{1} and
\begin{equation}\label{2}
\exp(-2\sqrt{M})\le\phi(z)\le \exp(2\sqrt{M}),\quad\forall\ z\in B_2.
\end{equation}
One way to verify that such solution $\phi$ exists is to define 
$$\phi(z)=\sup\{\varphi(z)\, :\, \Delta\varphi-V\varphi\ge 0\;\text{in}\; B_2,\; \varphi\le \exp(2\sqrt{M})\;\text{on}\;\partial B_2\}.$$
Then $\phi$ solves \eqref{1} and satisfies \eqref{2}. We can also see that $\phi$ is Lipschitz. Moreover, from the gradient estimate for Poisson's equation (see, for example, \cite{gt83}), we have that for $0<a_1<a_2$ with $a_2r<2$
\begin{equation}\label{interior}
\|\nabla\phi\|_{L^{\infty}(B_{a_1r})}\le\frac{CM}{r}\|\phi\|_{L^{\infty}(B_{a_2r})},
\end{equation}
where $C$ is an absolute constant.

If we set $u=\phi v$, then $v$ satisfies
\begin{equation}\label{3}
\nabla\cdot(\phi^2\nabla v)=0,\quad\text{in}\quad B_2.
\end{equation}
Let $\tilde v$ with $\tilde v(0)=0$ be the stream function related to $v$, i.e.,
\begin{equation}\label{3-1}
\left\{\begin{aligned}
\d_y\tilde v&=\phi^2\d_xv,\\ 
-\d_x\tilde v&=\phi^2\d_yv.
\end{aligned}\right.
\end{equation}
Let $g=\phi^2v+i\tilde v$, then $g$ satisfies
\begin{equation}\label{5}
\bar\d g=\bar\d\phi^2 v=\frac{\bar\d\phi^2}{2\phi^2}(g+\bar g)\quad\text{in}\quad B_2.
\end{equation}
As usual, we define $\bar\d=\frac{1}{2}(\d_x+i\d_y)$ and $\d=\frac 12(\d_x-i\d_y)$. Let
\begin{equation}\label{alpha}
\alpha=\frac{\bar\d\phi^2}{2\phi^2}=\frac{\bar\d\phi}{\phi}=\bar\d\log(\phi),
\end{equation}
then \eqref{5} is equivalent to
\begin{equation}\label{6}
\bar\d g=\alpha g+\alpha\bar g\quad\text{in}\quad B_2.
\end{equation}
We perform one more reduction. Defining
\begin{equation}\label{talpha}
\tilde \alpha=\begin{cases}
\alpha+\alpha{\bar g}/{g},\quad\text{if}\quad g\ne 0,\\
0,\quad\text{otherwise},
\end{cases}
\end{equation}
\eqref{6} now is reduced to
\begin{equation}\label{10}
\bar\partial g=\tilde\alpha g\quad\text{in}\quad B_2.
\end{equation}

We can solve for \eqref{10} directly. Before doing so, we need to obtain a precise estimate of $\alpha=\bar\d\log(\phi)$. In view of \eqref{2}, if we denote $\psi=\log\phi$, then $\psi$ satisfies
\begin{equation}\label{12}
|\psi(z)|\le 2\sqrt{M}\quad\text{in}\quad B_2,
\end{equation}
and solves the following equation
\begin{equation}\label{14}
\Delta\psi+|\nabla\psi|^2=V\quad\text{in}\quad B_2.
\end{equation}
The following estimate of $\nabla\psi$ is crucial.
\begin{lemma}\label{l1}
\begin{equation}\label{111}
\|\nabla\psi\|_{L^{\infty}(B_{7/5})}\le C\sqrt{M},
\end{equation}
where $C>0$ is an absolute constant.
\end{lemma}
\pf\ We begin with an $L^2$ estimate. Let $\theta\in C_0^{\infty}(B_2)$ satisfies $0\le\theta\le 1$ and $\theta=1$ for $(x,y)\in B_{9/5}$. Multiplying both sides of \eqref{14} by $\theta$, using \eqref{12}, and the integration by parts, we have that
\[
\int\theta|\nabla\psi|^2=\int \theta V-\int(\Delta\theta) \psi\le C(M+\sqrt{M})
\]
i.e.,
\begin{equation*}
\int_{B_{9/5}}|\nabla\psi|^2\le CM,
\end{equation*}
where $C>0$ is an absolute constant.

To proceed further, we rescale the equation \eqref{14}. Define $\vphi=\psi/C\sqrt{M}$ for some $C>0$. Then \eqref{14} becomes
\begin{equation}\label{114}
\eps\Delta\vphi+|\nabla\vphi|^2=\tilde V,
\end{equation}
where $\eps=1/C\sqrt{M}$ and $\tilde V=V/C^2M$. We can choose $C>10$ (so $\eps<1/10$) sufficiently large so that
\begin{equation}\label{1145}
\|\tilde V\|_{L^{\infty}(B_2)}\le 1,\quad\|\vphi\|_{L^{\infty}(B_2)}\le 1\quad\text{and}\quad\int_{B_{9/5}}|\nabla\vphi|^2\le 1.
\end{equation}
\begin{claim}\label{cl1}
For any $z\in B_{7/5}$ and $\eps<r<1/5$, we have
\[
\int_{B_r(z)}|\nabla\vphi|^2\le Cr^2.
\]
\end{claim}
\noindent{\bf Proof of Claim~\ref{cl1}}. It suffices to take $z=0$. Choose a cutoff function $\eta\in C_0^{\infty}(B_{2r})$ and $\eta=1$ on $B_r$. Denote $m=\int_{B_{2r}}\vphi/|B_{2r}|$. Clearly, we have
\begin{equation}\label{1401}
\begin{aligned}
0&=\eps\int\Delta((\vphi-m)\eta^2)=\eps\int\Delta\vphi\eta^2+4\eps\int\nabla\vphi\cdot\nabla\eta\eta+\eps\int\Delta\eta^2(\vphi-m)\\
&:=I+II+II.
\end{aligned}
\end{equation}
Now substituting $\eps\Delta\vphi$ from \eqref{114} in $I$ yields
\begin{equation}\label{1402}
I=-\int|\nabla\vphi|^2\eta^2+\int\tilde V\eta^2\le-\int|\nabla\vphi|^2\eta^2+Cr^2.
\end{equation}
Next we can estimate
\begin{equation}\label{1403}
|II|\le 4\eps\left(\int|\nabla\vphi|^2\eta^2\right)^{1/2}\left(\int|\nabla\eta|^2\right)^{1/2}\le\frac 12\int|\nabla\vphi|^2\eta^2+C\eps^2.
\end{equation}
Finally, for $III$, we obtain that
\begin{equation}\label{1405}
\begin{aligned}
|III|&\le C\eps r^{-2}\int_{B_{2r}}|\vphi-m|\le C\left(\int_{B_{2r}}|\vphi-m|^2\right)^{1/2}\left(\int_{B_{2r}}\eps^2r^{-4}\right)^{1/2}\\
&\le Cr\left(\int_{B_{2r}}|\nabla\vphi|^2\right)^{1/2}(\eps^2r^{-2})^{1/2}\le C\eps^2+\frac{1}{400}\int_{B_{2r}}|\nabla\vphi|^2,
\end{aligned}
\end{equation}
where we used the Poincar\'e inequality in the third inequality above. 

Putting \eqref{1401}-\eqref{1405} together  gives
\begin{equation}\label{1406}
\int_{B_r}|\nabla\vphi|^2\le C\eps^2+Cr^2+\frac{1}{200}\int_{B_{2r}}|\nabla\vphi|^2\le Cr^2+\frac{1}{200}\int_{B_{2r}}|\nabla\vphi|^2.
\end{equation}
Now if $r^2\ge1/200$, then \eqref{1406} implies
\[
\int_{B_r}|\nabla\vphi|^2\le Cr^2
\]
by the last estimate of \eqref{1145}. On the other hand, if $r^2<1/200$, we can pick a $k\in{\mathbb N}$ such that
\[
\frac 15\le 2^kr\le \frac 25,
\]
i.e., $r^2\ge(1/100)^k$. It is observed that $B_{2^kr}(z)\subset B_{9/5}$ for $z\in B_{7/5}$. Iterating \eqref{1406} in $k$ steps yields
\[
\int_{B_r}|\nabla\vphi|^2\le Cr^2+\frac{1}{200^k}\int_{B_{2^kr}}|\nabla\vphi|^2\le Cr^2.
\]
This ends of the proof of claim.\eproof

We will use Claim~\ref{cl1} to give a pointwise bound of $\nabla\vphi(z)$ for $z\in B_{7/5}$. For this end, we use another scaling. Let $\vphi_{\eps}(z)=\frac{1}{\eps}\vphi(\eps z)$, then we have that
\[
\nabla\vphi_{\eps}(z)=\nabla\vphi(\eps z),\quad\Delta\vphi_{\eps}(z)=\eps\Delta\vphi(\eps z)
\]
and therefore
\begin{equation}\label{1408}
\Delta\vphi_{\eps}+|\nabla\vphi_{\eps}|^2=\tilde V(\eps z):=\tilde V_{\eps}(z)\quad\text{for}\quad z\in B_2,
\end{equation}
\[
\|\tilde V_{\eps}\|_{L^{\infty}(B_2)}\le 1.
\]
Moreover, it follows from Claim~\ref{cl1} that
\[
\begin{aligned}
\int_{B_2}|\nabla\vphi_{\eps}|^2=\int_{B_2}|\nabla\vphi(\eps z)|^2=\frac{1}{\eps^2}\int_{B_{2\eps}}|\nabla\vphi|^2\le\frac{1}{\eps^2}\cdot C\eps^2=C.
\end{aligned}
\]
Now applying the elliptic regularity theorem to \eqref{1408} (see \cite[Chapter V, Theorem~2.3 and Proposition~2.1]{gi83}), we obtain that there exists $p>2$ such that
\[
\|\nabla\vphi_{\eps}\|_{L^p(B_1)}\le C.
\]
Define a new function
\[
\tilde{\vphi}_{\eps}(z)=\vphi_{\eps}(z)-\frac{1}{|B_1|}\int_{B_1}\vphi_{\eps}.
\]
Then $\nabla\vphi_{\eps}=\nabla\tilde{\vphi}_{\eps}$ and $\tilde{\vphi}_{\eps}$ satisfies
\[
\Delta\tilde{\vphi}_{\eps}=-|\nabla\tilde{\vphi}_{\eps}|^2+\tilde V_{\eps}:=\zeta\quad\text{in}\quad B_1.
\]
It is clear that $\|\zeta\|_{L^{p/2}(B_1)}\le C$. Moreover, from Poincar\'e's inequality, we have
\[
\|\tilde{\vphi}_{\eps}\|_{L^{p/2}(B_1)}\le\|\tilde{\vphi}_{\eps}\|_{L^p(B_1)}\le C\|\nabla\vphi_{\eps}\|_{L^p(B_1)}\le C.
\]
Elliptic regularity theorem implies
\[
\|\tilde{\vphi}_{\eps}\|_{W^{2,p/2}(B_r)}\le C
\]
for $r<1$. By bootstrapping arguments, we obtain that
\[
\|\nabla\tilde{\vphi}_{\eps}\|_{L^{\infty}(B_{r'})}=\|\nabla\vphi_{\eps}\|_{L^{\infty}(B_{r'})}=\|\nabla\vphi\|_{L^{\infty}(B_{\eps r'})}\le C
\]
with $r'<r$. The method clearly works for any $x\in B_{7/5}$. The derivation of \eqref{111} is then completed.\eproof
\begin{remark}\label{psirem}
In view of the definition of $\psi=\log\phi$, if we normalize $\phi$ at a point $\hat z\in B_{7/5}$ to be $1$, i.e., $\phi(\hat z)=1$, then the Lipschitz estimate of $\psi$ implies
\[
|\psi(z)|=|\psi(z)-\psi(\hat z)|\le C\sqrt{M}|z-\hat z|\quad\text{for}\quad z\in B_{7/5}.
\]
Therefore, we can see that if $\phi(\hat z)=1$, then
\[
\frac{1}{C}\le\phi(z)\le C,\quad\forall\ z\in B_{\frac{C}{\sqrt{M}}}(\hat z).
\]

\end{remark}

Using \eqref{111}, we have that $\|\alpha\|_{L^{\infty}(B_{7/5})}\le C\sqrt{M}$, which immediately implies
\begin{equation}\label{25}
\|\tilde\alpha\|_{L^{\infty}(B_{7/5})}\le C\sqrt{M}.
\end{equation}
Let $w(z)$, $z=x+iy$, be defined by
\[
w(z)=-\frac{1}{\pi}\int_{B_{7/5}}\frac{\tilde\alpha(\xi)}{\xi-z}d\xi,
\]
then $\bar\d w=\tilde\alpha$ in $B_{7/5}$ and
\begin{equation}\label{26}
\|w\|_{L^{\infty}(B_{7/5})}\le C\|\tilde\alpha\|_{L^{\infty}(B_{7/5})}\le C\sqrt{M}.
\end{equation}
It is clear that any solution $g$ of \eqref{10} in $B_{7/5}$ is given by
\[
g(z)=\exp(w(z))h(z),
\]
where $h$ is holomorphic in $B_{7/5}$.

Applying Hadamard's three-circle theorem to analytic function $h$, we have that
\begin{equation}\label{E120}
\|h\|_{L^{\infty}(B_{r_1})} \le (\|h\|_{L^{\infty}(B_{r/2})})^{\theta}(\|h\|_{L^{\infty}(B_{r_2})})^{1-\theta},
\end{equation}
where $r/2<r_1<r_2<7/5$ and
\begin{equation}\label{thetaexp}
\theta=\frac{\log(\frac{r_2}{r_1})}{\log(\frac{2r_2}{r})}.
\end{equation}
Substituting $h=\exp(-w)g$ into \eqref{E120} and using \eqref{26} implies
\begin{equation}\label{30}
\|g\|_{L^{\infty}(B_{r_1})} \le \exp(C\sqrt{M}) (\|g\|_{L^{\infty}(B_{r/2})})^{\theta}(\|g\|_{L^{\infty}(B_{r_2})})^{1-\theta}.
\end{equation}
Recall that $g=\phi^2 v+i\tilde v=\phi u+i\tilde v$ and hence, 
\begin{equation}\label{32}
|\phi u|\le|g|\le(|\phi u|+|\tilde v|).
\end{equation}
Following from $\tilde v(0)=0$ and \eqref{3-1} , for $z\in B_2$, we have that 
\begin{equation}\label{35}
|\tilde v(z)|\le Cr\exp(2\sqrt{M})\|\nabla v\|_{L^{\infty}(B_r)}\quad\forall\ z\in B_r,\ r<2.
\end{equation}
Assume that $u$ satisfies
\begin{equation}\label{11-2}
\|u\|_{L^{\infty}(B_2)}\le\exp(C_0\sqrt{M})
\end{equation}
for some $C_0>0$. Now choosing $r/2<1<6/5$, i.e., $r_1=1$, $r_2=6/5$ in \eqref{E120}, using \eqref{32}, $v=u/\phi$, and the interior estimate \eqref{interior} for $\phi$ and $u$, we conclude that
\begin{equation*}
\|u\|_{L^{\infty}(B_{1})}\le \exp(C\sqrt{M})(\|u\|_{L^{\infty}(B_r)})^{\theta}(\|u\|_{L^{\infty}(B_{7/5})})^{1-\theta}\le \exp(C\sqrt{M})(\|u\|_{L^{\infty}(B_r)})^{\theta},
\end{equation*}
where $C$ depends on $C_0$. In summary, we have proved that
\begin{theorem}\label{tt1}
Let $u$ be a real solution to \eqref{1}. Assume that \eqref{1-1}, \eqref{11-2} hold and furthermore
\[
\|u\|_{L^{\infty}(B_{1})}\ge 1.
\]
Then we have that
\[
\|u\|_{L^{\infty}(B_r)}\ge r^{C\sqrt{M}},
\]
where $C$ depends on $C_0$.
\end{theorem}

For any bounded solution $u$ solving
\begin{equation}\label{whosch}
\Delta u-Vu=0\quad\text{in}\quad\R^2,
\end{equation}
with $V\ge 0$ and $\|V\|_{L^{\infty}(\R^2)}\le 1$, Landis' conjecture (qualitative) is trivial. Nonetheless, using the scaling argument of \cite{bk05} as in the proof of Theorem~\ref{t0}, Theorem~\ref{tt1} immediately implies a quantitative version of  Landis's conjecture.
\begin{theorem}\label{tt2}
Let $u$ be a real solution to \eqref{whosch} with $V\ge 0$. Assume that $|u(z)|\le \exp(C_0|z|)$, $|V(z)|\le 1$, and $u(0)=1$. Then $u$ satisfies 
\begin{equation}\label{50}
\inf_{|z_0|=R}\sup_{|z-z_0|<1}|u(z)|\ge \exp(-CR\log R)\quad\text{for}\quad R\gg 1,
\end{equation}
where $C$ depends on $C_0$.
\end{theorem}

\section{Proofs of Theorem~\ref{t1}, \ref{t2}}\label{sec3}

Now we consider
\begin{equation}\label{2nd}
Lu:=\Delta u-\nabla(W(x,y)u)-V(x,y)u=0\quad\text{in}\quad B_2,
\end{equation}
where $W=(W_1,W_2)$ and $V$ are real-valued, measurable,  and $V(x,y)\ge 0$. As above, assume that
\[
\|V\|_{L^{\infty}(B_2)}\le M,\quad\|W\|_{L^{\infty}(B_2)}\le K
\]
with $M\ge 1$, $K\ge 1$. To construct a positive multiplier for \eqref{2nd}, we consider the adjoint operator of $L$, i.e.,
\begin{equation}\label{2-200}
L^{\ast}u=\Delta u+W\cdot\nabla u-Vu.
\end{equation}
Let $\phi_1=\exp({\lambda} x)$, then $\nabla\phi_1={\lambda}\exp({\lambda}x)(1,0)$, and
\[
\begin{aligned}
L^{\ast}\phi_1=&\Delta\phi_1+W\cdot\nabla\phi_1-V\phi_1\\
=&(\lambda^2-V)\exp({\lambda}x)+{\lambda}\exp({\lambda}x)W_1\\
\ge&(\lambda^2-\|V\|_{L^{\infty}(B_2)}-{\lambda}\|W\|_{L^{\infty}(B_2)})\phi_1.
\end{aligned}
\]
Hence, if $\lambda= (\sqrt{M}+K)$, then $L^{\ast}\phi_1\ge 0$. On the other hand, let $\phi_2=\exp(2(\sqrt{M}+K))$, then $L^{\ast}\phi_2=-V\phi_2\le 0$. Therefore, as before, there exists a positive Lipschitz solution $\phi$ satisfying $L^{\ast}\phi=0$  in $B_2$ and estimates
\begin{equation*}
\exp(-2(\sqrt{M}+K))\le\phi\le \exp(2(\sqrt{M}+K)),\quad\forall\ (x,y)\in B_2.
\end{equation*}

As above, if we let $v=u/\phi$, then $v$ satisfies
\begin{equation}\label{3-3}
\nabla\cdot(\phi^2(\nabla v-Wv))=0\quad\text{in}\quad B_2.
\end{equation}
Denote $\tilde v$ with $\tilde v(0)=0$ the stream function corresponding to $v$, i.e.
\begin{equation}\label{3-4}
\left\{\begin{aligned}
\d_y\tilde v&=\phi^2\d_xv-\phi^2W_1v,\\ 
-\d_x\tilde v&=\phi^2\d_yv-\phi^2W_2v.
\end{aligned}\right.
\end{equation}
As before, let $g=\phi ^2v+i\tilde v$, then $g$ solves
\begin{equation}\label{3-5}
\bar\d g=\gamma(g+\bar g)\quad\text{in}\quad B_2,
\end{equation}
where
\[
\gamma=\frac{\bar\d\phi^2}{2\phi^2}+\frac 12(W_1+iW_2)=\bar\d\log\phi+\frac 12(W_1+iW_2).
\]
Likewise, let
\[
\tilde \gamma=\begin{cases}
\gamma+\gamma{\bar g}/{g},\quad\text{if}\quad g\ne 0,\\
0,\quad\text{otherwise},
\end{cases}
\]
then \eqref{3-5} becomes
\begin{equation}\label{3-10}
\bar\partial g=\tilde\gamma g\quad\text{in}\quad B_2.
\end{equation}

Let $\phi=\exp(\psi)$, then $\psi$ satisfies
\begin{equation}\label{0417}
|\psi(z)|\le 2(\sqrt{M}+K),\;\;\text{for all}\;\; z\in B_2
\end{equation}
and
\begin{equation}\label{3-12}
\Delta\psi+|\nabla\psi|^2+W\cdot\nabla\psi=V\quad\text{in}\quad B_2.
\end{equation}
Similarly, we prove the following estimate of $\nabla\psi$.
\begin{lemma}\label{l3-1}
\begin{equation}\label{3-15}
\|\nabla\psi\|_{L^{\infty}(B_{7/5})}\le C(\sqrt{M}+K).
\end{equation}
\end{lemma}
\pf This lemma can be proved in the same way as in Lemma~\ref{l1}. We first derive an $L^2$ bound. Let $\theta\in C_0^{\infty}(B_2)$ satisfies $0\le\theta\le 1$ and $\theta=1$ for $(x,y)\in B_{9/5}$. Multiplying both sides of \eqref{3-12} by $\theta$, using \eqref{0417} and the integration by parts, we obtain that
\[
\int\theta|\nabla\psi|^2=\int V\theta-\int \theta W\cdot\nabla\psi-\int\Delta\theta\psi\le\frac 12\int\theta|\nabla\psi|^2+C(M+K^2+\sqrt{M}+K)
\]
i.e.,
\[
\int_{B_{9/5}}|\nabla\psi|^2\le C(M+K^2).
\]
To bound $\nabla\psi(x)$ for all $x\in B_{7/5}$ by $C(\sqrt{M}+K)$, choosing $\eps=1/C(\sqrt{M}+K)$ and $\vphi=\psi/C(\sqrt{M}+K)$, we proceed as in the proof of Lemma~\ref{l1}. We avoid repeating the arguments here.\eproof

With the help of estimate \eqref{3-15}, we have that
\[
\|\tilde\gamma\|_{L^{\infty}(B_{7/5})}\le C(\sqrt{M}+K).
\]
Therefore, any solution $g$ of \eqref{3-10} in $B_{7/5}$ is represented by
\[
g=\exp(\tilde w)h,
\]
where $h$ is holomorphic in $B_{7/5}$ and 
\[
\|\tilde w\|_{L^{\infty}(B_{7/5})}\le C(\sqrt{M}+K).
\]
The remaining arguments in proving Theorem~\ref{t1} and \ref{t2} are exactly similar to those of Theorem~\ref{tt1} and \ref{tt2}.

\section{Proofs of Theorem~\ref{t3}, \ref{t4}}\label{sec4}

Recall that $u$ is a real solution of
\begin{equation}\label{5-8}
\Delta u+W\cdot\nabla u-Vu=0\quad\text{in}\quad B_2,
\end{equation}
where $W=(W_1,W_2)$ and $V\ge 0$ are real-valued and
\[
\|W\|_{L^{\infty}(B_2)}\le K\quad\text{and}\quad\|V\|_{L^{\infty}(B_2)}\le M.
\]
Let $\phi$ be the positive solution of \eqref{5-8} constructed in Section~\ref{sec3} for $L^{\ast}u=0$ (see \eqref{2-200}). Defining $v=u/\phi$ implies that $v$ satisfies
\begin{equation}\label{5-10}
\Delta v+(2\nabla\psi+W)\cdot\nabla v=0\quad\text{in}\quad B_2.
\end{equation}
In view of estimate \eqref{3-15}, we have 
\[
\|2\nabla\psi+W\|_{L^{\infty}(B_{7/5})}\le C(\sqrt{M}+K).
\]
Note that $4\Delta=\bd\d$. Hence, \eqref{5-10} can be written as
\begin{equation}\label{55-10}
\bd(\d v)=\widetilde W(\d v)\quad\text{in}\quad B_2
\end{equation}
with
\[
\|\widetilde W\|_{L^{\infty}(B_{7/5})}\le C(\sqrt{M}+K).
\]
Argued as above, using Hadamard's three-circle theorem, we have that
\begin{equation}\label{5-12}
\|\nabla v\|_{L^{\infty}(B_{r_1})} \le \exp(C(\sqrt{M}+K))(\|\nabla v\|_{L^{\infty}(B_{r/2})})^{\theta}(\|\nabla v\|_{L^{\infty}(B_{r_2})})^{1-\theta},
\end{equation}
where $r/2<r_1$, $r_1=6/5$, $r_2=7/5$, and
\[
\theta=\frac{\log(\frac{r_2}{r_1})}{\log(\frac{2r_2}{r})}.
\]
Using the interior estimate again, \eqref{5-12} becomes
\begin{equation}\label{5-121}
\|\nabla v\|_{L^{\infty}(B_{6/5})} \le \exp(C(\sqrt{M}+K))(r^{-1}\|u\|_{L^{\infty}(B_{r})})^{\theta}.
\end{equation}

Now we would like to bound the left hand side of \eqref{5-121} from below using the a priori condition $\|u\|_{L^{\infty}(B_1)}\ge 1$. Since $\|u\|_{L^{\infty}(B_1)}\ge 1$, there exists $z_0\in B_1$ such that $|u(z_0)|\ge 1$. It suffices to assume $u(z_0)\ge1$. For a real-valued $u$, it is clear that given any $a>0$, either $u(z)\ge a$ for all $z\in B_{6/5}$ or there exists $z_1\in B_{6/5}$ such that $u(z_1)<a$. Here we would like to choose an appropriate $a$. For the latter case, recalling that
\[
\exp(-2(\sqrt{M}+K))\le\phi(z)\le\exp(2(\sqrt{M}+K)),
\]
we can see that
\[
\frac{u(z_1)}{\phi(z_1)}\le\frac{a}{\phi(z_1)}\le a\exp(2(\sqrt{M}+K))
\]
and
\[
\frac{u(z_0)}{\phi(z_0)}\ge\exp(-2(\sqrt{M}+K)).
\]
Therefore, if we set
\[
a\exp(2(\sqrt{M}+K))=\frac{1}{2}\exp(-2(\sqrt{M}+K)),\;\;\text{i.e.},\;\; a=\frac 12\exp(-4(\sqrt{M}+K)),
\]
then we have that
\[
\frac{u(z_1)}{\phi(z_1)}\le\frac 12\exp(-2(\sqrt{M}+K))
\]
and thus
\[
\|\nabla v\|_{L^{\infty}(B_{6/5})}\ge|v(z_0)-v(z_1)|\ge\frac{u(z_0)}{\phi(z_0)}-\frac{u(z_1)}{\phi(z_1)}\ge\frac 12\exp(-2(\sqrt{M}+K)),
\]
which implies \eqref{r18} with the help of \eqref{5-121}. Now for the other case, i.e., $u(z)\ge a$ for all $z\in B_{6/5}$, \eqref{r18} is obviously satisfied. This completes the proof of Theorem~\ref{t3}. By Theorem~\ref{t3}, the proof of Theorem~\ref{t4} follows from the usual scaling augment in \cite{bk05}.

\section{Landis' conjecture in an exterior domain}\label{sec5}

In this section we prove Landis' conjecture in an exterior domain, Theorem~\ref{t5}. Recall that we consider
\begin{equation}\label{equ0}
\Delta u-V(x,y)u=0\quad\text{in}\quad B^c_1.
\end{equation}
Assume that the potential $V$ is defined everywhere in $\R^2$ and satisfies $V\ge 0$ and 
\begin{equation}\label{vpo}
\|V\|_{L^{\infty}(\R^2)}\le 1.
\end{equation}
As mentioned in the Introduction, for \eqref{equ0}, the local vanishing order result, Theorem~\ref{tt1}, remains true. However, the scaling argument fails to imply Landis' conjecture. We have to work harder to prove the conjecture in this case. Here our main tool is a Carleman estimate. But we need to set up everything carefully before applying the Carleman estimate

Let $z_0'\in\R^2$ with $|z_0'|\gg 1$. Since \eqref{equ0} is invariant under rotation, we can assume that $z'_0=|z_0'|e_1$, where $e_1=(1,0)$. Translating the origin to $-5e_1/2$, \eqref{equ0} becomes
\begin{equation}\label{equ1}
\Delta u-V(x,y)u=0\quad\text{in}\quad B^c_1(-5e_1/2).
\end{equation}
By abuse of notation, we continue to write $u$ and $V$ in the equation in the new coordinates. Now we denote $z_0=(|z_0'|-5/2)e_1$ and set $R=|z_0|$.  As before, we define the scaled solution $u_R(z)=u(ARz+z_0)$, where $A>0$ will be determined later. Therefore, $u_R$ solves
\begin{equation}\label{equ2}
\Delta u_R-V_Ru_R=0\quad\text{in}\quad B_{\frac{1}{AR}}^c(z_1),
\end{equation}
where
\[
z_1=-(\frac 1A+\frac{5}{2AR})e_1
\]
and $V_R(z)=(AR)^2V(ARz+z_0)$, thus,
\[
\|V_R\|_{L^{\infty}(\R^2)}\le (AR)^2.
\]
Under the above change of coordinates, the origin moves to
\[
\hat z=-\frac{z_0}{AR}=-\frac 1Ae_1.
\] 
We choose a large $A$ so that
\[
B_{\frac{1}{AR}}(z_1)\subset B_{7/5}.
\]

Let $\phi$ be the positive solution of 
\begin{equation}\label{phihat}
\Delta\phi-V_R\phi=0\quad\text{in}\quad B_2
\end{equation}
as constructed in Section~2. Here we normalize $\phi$ at $\hat z$, i.e., $\phi(\hat z)=1$. Likewise, we define $\psi=\log\phi$. Then $\psi$ satisfies
\begin{equation}\label{esti1}
\|\psi\|_{L^{\infty}(B_{7/5})}\le C(AR)\quad\text{and}\quad\|\nabla\psi\|_{L^{\infty}(B_{7/5})}\le C(AR).
\end{equation}
To simplify the notation, we suppress the subscript $R$ of $u_R$ and denote $u_R$ by $u$. As before, let $v=u/\phi$, then $v$ solves
\begin{equation}\label{equ3}
\nabla\cdot(\phi^2\nabla v)=0\quad\text{in}\quad B_2\setminus B_{\frac{1}{AR}}(z_1).
\end{equation}
Note that here the domain of \eqref{equ3} is not simply-connected. The stream function of $v$, solution to \eqref{equ3}, may not exist. However, we can construct an "approximate" stream function of $v$. To do this, we choose a cutoff function $\chi\equiv 1$ on $|z-z_1|\ge \frac{9}{8AR}$ and $\equiv 0$ for $|z-z_1|\le\frac{17}{16AR}$. Note that $\nabla\chi$  is supported on $\frac{17}{16AR}\le|z-z_1|\le\frac{9}{8AR}$. Next we denote
\[
ae_1=\left(\frac{9}{8AR}-(\frac 1A+\frac{5}{2AR})\right)e_1=\left(-\frac 1A-\frac{11}{8AR}\right)e_1,
\]
i.e.,
\[
a=-\frac 1A-\frac{11}{8AR}.
\]
Now we define 
\[
\tilde v(x,y)=\int_{a}^x-[\chi\phi^2\d_yv](s,y)ds+\int_0^y[\chi\phi^2\d_xv](a,s)ds.
\]
It is easy to see that
\begin{equation}\label{equ5}
\begin{aligned}
\d_y\tilde v(x,y)&=(\chi\phi^2\d_xv)(a,y)+\int_{a}^x-[\d_y\chi\phi^2\d_yv](s,y)ds+\int_{a}^x[\chi\d_x(\phi^2\d_xv)](s,y)ds\\
&=\chi\phi^2\d_xv(x,y)+\int_{a}^x-[\d_y\chi\phi^2\d_yv+\d_x\chi\phi^2\d_xv](s,y)ds
\end{aligned}
\end{equation}
and
\begin{equation}\label{equ6}
\d_x\tilde v(x,y)=-\chi\phi^2\d_yv(x,y).
\end{equation}
As before, we set $g=\chi \phi^2v+i\tilde v$ and hence
\begin{equation}\label{equ8}
\begin{aligned}
\bd g&=(\bd\phi^2)\chi v+(\bd\chi)\phi^2v+\frac{1}{2}\int_{a}^x[\d_y\chi\phi^2\d_yv+\d_x\chi\phi^2\d_xv](s,y)ds,\\
&=\frac{\bd\phi^2}{2\phi^2}(g+\bar g)+(\bd\chi)\phi u+\frac{1}{2}\int_{a}^x[\d_y\chi\phi^2\d_yv+\d_x\chi\phi^2\d_xv](s,y)ds\\
&=\alpha g+\alpha\bar g+(\bd\chi)\phi u+\frac{1}{2}\int_{a}^x[\d_y\chi\phi^2\d_yv+\d_x\chi\phi^2\d_xv](s,y)ds,
\end{aligned}
\end{equation}
where $\alpha$ is given in \eqref{alpha}. Defining $\tilde\alpha$ as in \eqref{talpha}, \eqref{equ8} now is equivalent to
\begin{equation}\label{equ9}
\bar\partial g=\tilde\alpha g+(\bd\chi)\phi u+\frac{1}{2}\int_{a}^x[\d_y\chi\phi^2\d_yv+\d_x\chi\phi^2\d_xv](s,y)ds\quad\text{in}\quad B_2.
\end{equation}

We now write $\hat z$ as a point in the complex plane, i.e., $\hat z=-\frac{1}{A}+i0$. Let $w(z)$ be defined by
\[
w(z)=\frac{1}{\pi}\int_{B_{7/5}}\frac{\tilde\alpha}{\xi-z}d\xi-\frac{1}{\pi}\int_{B_{7/5}}\frac{\tilde\alpha}{\xi-\hat z}d\xi,
\]
then $\bar\d w=-\tilde\alpha$ in $B_{7/5}$. Recall that
\[
\|\tilde\alpha\|_{L^{\infty}(B_{7/5})}\le C(AR).
\]
In view of \cite[(6.9a)]{ve62}, we have the following estimate of $w(z)$.
\begin{equation}\label{wze}
|w(z)|\le C(AR)|z-\hat z|\log\left(\frac{C}{|z-\hat z|}\right),\quad\forall\ z\in B_{7/5}.
\end{equation}
Let $h=e^wg$, then it follows from \eqref{equ9} that
\begin{equation}\label{equ11}
\begin{aligned}
\bd h&=e^{w}(\bd\chi)\phi u+\frac{e^{w}}{2}\int_{a}^x[\d_y\chi\phi^2\d_yv+\d_x\chi\phi^2\d_xv](s,y)ds\\
&=H_1+H_2\quad\text{in}\quad B_{7/5},
\end{aligned}
\end{equation}
where
\[
H_1=e^{w}(\bd\chi)\phi u\;\;\text{and}\;\; H_2=\frac{e^{w}}{2}\int_{a}^x[\d_y\chi\phi^2\d_yv+\d_x\chi\phi^2\d_xv](s,y)ds.
\]

We now come to the Carleman estimate. Here we will use the following estimate for $\bd$ from \cite[Proposition~2.1]{df90}. Let $\vphi_{\tau}(z)=\vphi_{\tau}(|z|)=-\tau\log|z|+|z|^2$, then for any $h\in C_0^{\infty}(B_{7/5}\setminus\{0\})$, we have that
\begin{equation}\label{care}
\int |\bd h|^2e^{\vphi_{\tau}}\ge\frac 14\int(\Delta\vphi_{\tau})|h|^2e^{\vphi_{\tau}}=\int|h|^2e^{\vphi_{\tau}}.
\end{equation}
Note that $\vphi_{\tau}$ is decreasing in $|z|$ for $\tau>8$. We introduce another cutoff function $0\le\zeta\le 1$ satisfying
\[
\zeta(z)=\left\{
\begin{aligned}
0,&\quad\text{when}\ |z|<\frac{1}{4AR},\\
1,&\quad\text{when}\ \frac{1}{2AR}<|z|<1,\\
0,&\quad\text{when}\ |z|>6/5.
\end{aligned}\right.
\] 
Hence the following estimates holds
\begin{equation}\label{estz}
|\nabla\zeta(z)|\le C(AR)\;\;\text{for}\;\; z\in X\;\;\text{and}\;\;|\nabla\zeta(z)|\le C\;\;\text{for}\;\; z\in Y,
\end{equation}
where
\[
X=\{\frac{1}{4AR}<|z|<\frac{1}{2AR}\}\;\; \text{and}\;\; Y=\{1<|z|<6/5\}.
\]
We also denote
\[
Z=\{\frac{1}{2AR}<|z|<1\}.
\]
The relative positions of different domains in the rescaled problem are shown in the following figure, Figure~\ref{fig1}.
\begin{figure}[ht]
\centering
\begin{tikzpicture}
\draw  (6,0) -- (-6,0);
\draw (5.2,0) circle [radius=.65];
\draw (.28,0) circle [radius=.65];
\draw[fill] (.28,0) circle [radius=0.04] node[below]{$\scriptscriptstyle -\frac{1}{A}$} node[above]{$\hat z$};
\draw[fill=blue] (-1.5,0) circle [radius=.65];
\draw[fill] (-1.5,0) circle [radius=0.04] node[below]{$\scriptscriptstyle -\frac{1}{A}-\frac{5}{2AR}$} node[above]{$z_1$};
\draw[fill] (-0.5,0) circle [radius=0.04]  node[above]{$a$};
\draw[line width=30pt, color=green] (-5.0,0) arc [radius=9.7, start angle=180, end angle=150];
\draw[line width=30pt, color=green] (-5.0,0) arc [radius=9.7, start angle=180, end angle=210];
\draw[fill] (-4.44,0) circle [radius=0.04] node[below right] {$-1$};
\draw[fill] (-5.53,0) circle [radius=0.04] node[below left] {$-\frac{6}{5}$};
\draw[line width=7pt, color=red] (-1.5,0) circle [radius=0.85];
\draw[line width=1pt, dotted] (.28,0) circle [radius=2.75];
\draw[line width=5pt, color=green] (5.2,0) circle [radius=0.4];
\draw[fill] (5.2,0) circle [radius=0.04] node[below]{$\scriptstyle 0$} node[above]{$z_0$};
\node at (.28,1.5) {size of ball $=\frac{C}{AR}$};
\end{tikzpicture}
\caption{The figure represents the domain in the rescaled problem. $\nabla\chi$ is supported in the red region and $\nabla\zeta$ is supported in the green region. All three balls centered at $z_0$, $\hat z$, and $z_1$ have radius $\frac{1}{AR}$, corresponding to balls of radius $1$ in the unscaled problem. $\phi$ is bounded from above and below (away from zero) with bounds  independent of $R$ in the ball of size $\frac{C}{AR}$ centered at $\hat z$.}
\label{fig1}
\end{figure}
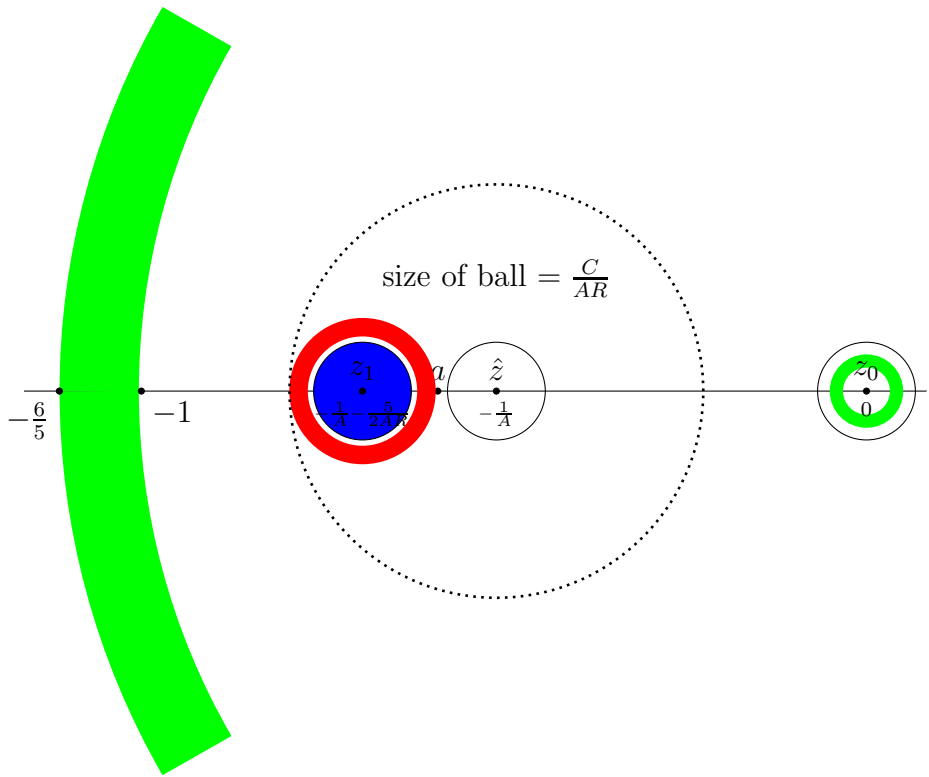

Applying the Carleman estimate \eqref{care} to $\zeta h$ gives
\begin{equation}\label{der1}
\begin{aligned}
\int_Z|h|^2e^{\vphi_{\tau}}&\le 2\int(|\bd\zeta h|^2+|\zeta\bd h|^2)e^{\vphi_{\tau}}\\
&\le C(AR)^2\int_X|h|^2e^{\vphi_{\tau}}+C\int_Y|h|^2e^{\vphi_{\tau}}+\int_{\widetilde Z}(|H_1|^2+|H_2|^2)e^{\vphi_{\tau}},
\end{aligned}
\end{equation}
where
\[
\widetilde Z=\{\frac{1}{4AR}<|z|<\frac 65\}.
\]
The terms $\int_Z|h|^2e^{\vphi_{\tau}}$ and $\int_{\widetilde Z}(|H_1|^2+|H_2|^2)e^{\vphi_{\tau}}$ are the most crucial ones in \eqref{der1}. We would like to study them in more detail. We begin with the first one. It is clear that for $A$ large, $R$ large
\[
\int_Z|h|^2e^{\vphi_{\tau}}\ge\int_{B_{\frac{1}{AR}}(\hat z)}|h|^2e^{\vphi_{\tau}}.
\]
From Remark~\ref{psirem}, we have that $0<C_1\le\phi(z)\le C_2$ in the ball $B_{\frac{1}{AR}}(\hat z)$ with absolute constants $C_1$, $C_2$. Since $\chi\equiv 1$ on $B_{1/AR}(\hat z)$, we thus obtain that
\begin{equation}\label{gcv}
|g|\ge |\phi^2 v|=|\phi u|\ge C|u|\;\;\text{for}\;\; z\in B_{\frac{1}{AR}}(\hat z).
\end{equation}
On the other hand, \eqref{wze} implies
\[
|w(z)|\le C(AR)\frac{1}{AR}\log({C}{AR})=C\log({C}{AR})\;\;\text{for}\;\; z\in B_{\frac{1}{AR}}(\hat z),
\]
i.e.,
\begin{equation}\label{ewz}
e^{w(z)}\ge\frac{1}{C(AR)^C}\;\;\text{for}\;\; z\in B_{\frac{1}{AR}}(\hat z).
\end{equation}
Combining \eqref{gcv}, \eqref{ewz}, and using that for $z\in B_{1/AR}(\hat z)$, $|z|\le\frac{1}{AR}+\frac 1A$,  we have
\begin{equation}\label{bar}
\int_{Z}|h|^2e^{\vphi_{\tau}}\ge\frac{e^{\vphi_{\tau}(\frac 1A+\frac{1}{AR})}}{C(AR)^C}\int_{B_{\frac{1}{AR}}(\hat z)}|u|^2.
\end{equation}

Next we look at $\int_{\widetilde Z}|H_1|^2e^{\vphi_{\tau}}$. Recall that $H_1$ is supported in $\frac{17}{16AR}\le|z-z_1|\le\frac{9}{8AR}$ since $\nabla\chi$ is. It is clear that
\[
G:=\{\frac{17}{16AR}\le|z-z_1|\le\frac{9}{8AR}\}\subset\{|z-\hat z|\le\frac{29}{8AR}\}.
\]
As in \eqref{ewz}, we can see that
\begin{equation}\label{ewz1}
e^{w(z)}\le{C(AR)^C}\;\;\text{for}\;\; \frac{17}{16AR}\le|z-z_1|\le\frac{9}{8AR}.
\end{equation}
Using \eqref{ewz1} and the same argument based on Remark~\ref{psirem} as above, we have that
\begin{equation}\label{h1e}
\int_{\widetilde Z}|H_1|^2e^{\vphi_{\tau}}\le C(AR)^C\int_{\frac{17}{16AR}\le|z-z_1|\le\frac{9}{8AR}}|u|^2e^{\vphi_{\tau}}\le C(AR)^Ce^{\vphi_{\tau}(\frac 1A+\frac{29}{8AR})}.
\end{equation}

Finally, we study $\int_{\widetilde Z}|H_2|^2e^{\vphi_{\tau}}$. Observe that
\[
\supp H_2\subset\{-\frac 65<x<-\frac 1A-\frac{11}{8AR},\ |y|<\frac{9}{8AR}\}.
\]
Argued as above, $\phi$ is uniformly bounded from above and below in $G$. In view of the a priori boundedness assumption on $u$, Lemma~\ref{l1}, and the interior estimate, we obtain that
\begin{equation}\label{h2e}
|H_2(z)|\le C(AR)^Ce^{|w(z)|}\quad\text{for}\quad z\in\supp H_2.
\end{equation}
For $z\in\supp H_2$, we have that
\[
|z-\hat z|\le |y|+|x+\frac{1}{A}|\le\frac{9}{8AR}+|x+\frac{1}{A}|,
\]
\[
\frac{11}{8AR}\le|x+\frac 1A|\le |z-\hat z|,
\]
and so \eqref{wze} implies
\begin{equation}\label{h2ee}
|H_2|\le\exp(C(AR)|x+\frac 1A|\log(AR)).
\end{equation}
We will multiply by $\exp({-\vphi_{\tau}(\frac 1A+\frac{1}{AR})})$ on both sides of \eqref{der1}. Thus, we want to take a closer look at
\[
\log\left(\frac{\frac 1A+\frac{1}{AR}}{|z|}\right)\quad\text{for}\quad z\in\supp H_2.
\]
Since 
\[
\frac{\frac 1A+\frac{1}{AR}}{|z|}<1,
\]
because $|x|\ge\frac{1}{A}+\frac{11}{8AR}$, it suffices to estimate
\[
1-\frac{1+\frac 1R}{|Az|}=\frac{A|z|-1-\frac 1R}{A|z|}\ge\frac{A|x|-1-\frac 1R}{A|z|}\ge \frac{A( |x|-\frac 1A)}{A|z|}\ge C(|x|-\frac 1A)
\]
for $z\in\supp H_2$. Here $C$ depends on $A$, but $A$ has been chosen. Note that when $z=(x,y)\in\supp H_2$, $|x+\frac 1A|=|x|-\frac 1A$. Thus, we have 
\[
\log\left(\frac{\frac 1A+\frac{1}{AR}}{|z|}\right)\le -C|x+\frac 1A|
\]
and 
\begin{equation}\label{h22e}
\exp(\vphi_{\tau}(z))\exp({-\vphi_{\tau}(\frac 1A+\frac{1}{AR})})|H_2|\le\exp(C(AR)|x+\frac 1A|\log(AR)-C\tau|x+\frac 1A|).
\end{equation}

Now we can put everything together. Multiplying $C(AR)^C\exp(-\vphi_{\tau}(\frac 1A+\frac{1}{AR}))$ on both sides of \eqref{der1}, using \eqref{bar}, \eqref{h1e}, \eqref{h22e}, the second inequality of \eqref{32}, interior estimates, and the a priori bound of $u$, we obtain that
\begin{equation}\label{h12e}
\begin{aligned}
\int_{B_{\frac{1}{AR}(\hat z)}}|u|^2&\le C(AR)^C\exp(C(AR)\log(AR))\frac{\exp(\vphi_{\tau}(\frac{1}{4AR}))}{\exp(\vphi_{\tau}(\frac 1A+\frac{1}{AR}))}\int_{B_{\frac{1}{AR}(0)}}|u|^2\\
&\;\;\;+C(AR)^C\exp(C(AR)\log(AR))\frac{\exp(\vphi_{\tau}(1))}{\exp(\vphi_{\tau}(\frac 1A+\frac{1}{AR}))}\\
&\;\;\;+C(AR)^C\frac{\exp({\vphi_{\tau}(\frac 1A+\frac{29}{8AR})})}{\exp(\vphi_{\tau}(\frac 1A+\frac{1}{AR}))}\\
&\;\;\;+\exp\left(\left(C(AR)\log(AR)-C\tau\right)\frac{11}{8AR}\right),
\end{aligned}
\end{equation}
where the last term of \eqref{h12e} was derived by using $|x+\frac 1A|\ge\frac{11}{8AR}$ and taking $\tau$ large enough such that
\begin{equation}\label{ttau}
C(AR)\log(AR)-C\tau<0.
\end{equation}
To fulfill \eqref{ttau}, it suffices to choose
\begin{equation}\label{tauu}
\tau=\tilde C(AR)\log(AR)
\end{equation}
for an appropriate large fixed constant $\tilde C$. Consequently, the last term of \eqref{h12e} satisfies
\begin{equation}\label{3term}
\exp\left(\left(C(AR)\log(AR)-C\tau\right)\frac{11}{8AR}\right)\le\exp(-\tilde C\log(AR))=(AR)^{-\tilde C}.
\end{equation}

Rescaling back to the original variables, we observe that
\[
\int_{B_{\frac{1}{AR}(\hat z)}}|u_R|^2=\frac{1}{(AR)^2}\int_{B_1(0)}|u|^2\ge \frac{C_0}{(AR)^2}\;\;(\text{from}\;\;\eqref{apriorib})\quad\text{and}\quad\int_{B_{\frac{1}{AR}(0)}}|u_R|^2=\frac{1}{(AR)^2}\int_{B_1(x_0)}|u|^2
\]
Finally, choosing $\tau$ as in \eqref{tauu} (choosing $\tilde C$ larger if necessary) and taking $R$ sufficiently large, it is not hard to see that
\[
\left\{\begin{aligned}
&C(AR)^{C+2}\exp(C(AR)\log(AR))\frac{\exp(\vphi_{\tau}(\frac{1}{4AR}))}{\exp(\vphi_{\tau}(\frac 1A+\frac{1}{AR}))}\le\exp(C(AR)(\log(AR))^2),\\
&C(AR)^{C+2}\exp(C(AR)\log(AR))\frac{\exp(\vphi_{\tau}(1))}{\exp(\vphi_{\tau}(\frac 1A+\frac{1}{AR}))}\to 0,\\
&C(AR)^{C+2}\frac{\exp({\vphi_{\tau}(\frac 1A+\frac{29}{8AR})})}{\exp(\vphi_{\tau}(\frac 1A+\frac{1}{AR}))}\to 0,\\
&(AR)^2\exp\left(\left(C(AR)\log(AR)-C\tau\right)\frac{11}{8AR}\right)\to 0\;\; (\text{from}\; \eqref{3term}).
\end{aligned}\right.
\]
Therefore, if $R$ is large enough, then the last three terms on the right hand side of \eqref{h12e} can be absorbed by the term on the left. The proof now is completed. 

The same method can be used to generalize Theorems~\ref{t2}, \ref{t3}, \ref{t4}, to the case of an exterior domain. Before ending this section, we would like to remark that we can express \eqref{equ3} as
\begin{equation}\label{equ66}
\Delta v+2\nabla\psi\cdot\nabla v=0\quad\text{in}\quad B_2\setminus B_{\frac{1}{AR}}(z_1),
\end{equation}
which is exactly in the same form as in \eqref{5-10}. Thus, using $g=\d v$, \eqref{equ66} can be transformed into a $\bd$ equation for $g$ in $B_2\setminus B_{\frac{1}{AR}}(z_1)$ without introducing an approximate stream function as before. However, this reduction does not work for the equation \eqref{e2} (see \eqref{3-3}). Our point here is to introduce a unified approach using an approximate stream function to prove Landis' conjecture (Theorems~\ref{t2}, \ref{t3}, \ref{t4}) in an exterior domain. 

\section{Conclusions and discussions}\label{conl}

We proved positive answers to Landis' conjecture for $\Delta u-\nabla(Wu)-Vu=0$ and $\Delta u+W\nabla u-Vu=0$ in the plane or in an exterior domain with or without $W$ under the assumption of $V\ge 0$. One key ingredient of our method is the existence of a global positive multiplier which convert the original equation into an equation in divergence form. Another important tool is a global three-ball inequality with an \emph{optimal} exponent ($\theta$ in \eqref{thetaexp}), namely, Hadamard's three-circle theorem. In the case of an exterior domain, we use a global Carleman estimate for $\bd$. For the general potential $V$, one is also able to construct a positive multiplier that can be used to convert the original equation into a divergence-form equation (see, for example, \cite{al98}, \cite{al10}, \cite{sc98}). However, this positive multiplier exists only in small balls, with radius depending on $\|V\|_{L^{\infty}}$. In other words, the divergence-form equation is valid only locally, which leads to a local quantitative uniqueness estimate (three-ball inequality). If we try to prove bounds like \eqref{r1} or \eqref{r2}, we need to iterate that local estimate, which will give rise to double exponential bounds and are worse than the results derived by the Carleman method.  

One trick to remove the assumption $V\ge 0$ is to consider the new function $u_{\xi}(x,y)=\cosh(\lambda\xi)u(x,y)$. Then if $u$ satisfies $\Delta_{x,y} u-Vu=0$, $u_{\xi}$ satisfies 
\begin{equation}\label{4-1}
\d_{\xi}^2u_{\xi}+\Delta_{x,y}u_{\xi}-(\lambda^2-V)u_{\xi}=0\quad\text{in}\quad\R^3.
\end{equation}
The new potential function $\lambda^2-V$ is non-negative if $\lambda\ge\sqrt{\|V\|_{L^{\infty}}}$. Consequently, we can construct a global positive solution to \eqref{4-1}, still denoted by $\phi$, satisfying similar estimates to the ones described in Section~\ref{sec2}. Using $\phi$, \eqref{4-1} becomes
\begin{equation}\label{4-2}
\nabla_{\xi,x,y}\cdot(\phi^2\nabla_{\xi,x,y}v_{\xi})=0\quad\text{in}\quad\R^3
\end{equation}
with $v_{\xi}=u_{\xi}/\phi$. However, at this stage, an optimal three-ball inequality or an optimal Carleman estimate with a suitable weight function for \eqref{4-2} are not available. Therefore, new ideas seem to be needed to resolve Landis' conjecture in the general case.


\begin{thebibliography}{99999}
\bibitem[Al98]{al98}
G. Alessandrini, \emph{On Courant's nodal domain theorem}, Forum Mathematicum, \textbf{10} (1998), 521-532. 

\bibitem[Al10]{al10}
G. Alessandrini, \emph{Strong unique continuation for general elliptic equations in 2D},  J. Math. Anal. App., \textbf{386} (2012), 669-676.

\bibitem[An58]{an58}
P. Anderson, \emph{Absence of diffusion in certain random lattices}, Phys. Review, \textbf{109} (1958), 1492-1505.


\bibitem[BK05]{bk05}
J. Bourgain and C. Kenig, \emph{On localization in the Anderson-Bernoulli model in higher dimensions}, Invent. Math., \textbf{161} (2005), 389-426.


\bibitem[DF88]{df88}
H. Donnelly and C. Fefferman, \emph{Nodal sets of eigenfunctions on
Riemannian manifolds}, Invent. Math. {\bf 93} (1988), 161-183.

\bibitem[DF90]{df90}
H. Donnelly and C. Fefferman, \emph{Nodal sets for eigenfunctions of the Laplacian on surfaces}, J. of AMS {\bf 3} (1990), 333-353.

\bibitem[Da12]{da12}
B. Davey, \emph{Some quantitative unique continuation results for eigenfunctions of the magnetic Schr\"odinger operator}, Comm. in PDE, {\bf 39} (2014), 876-945.


\bibitem[GL86]{gl86}
N. Garofalo and F.H. Lin, \emph{Monotonicity properties of
variational integrals, $A_p$ weights and unique continuation},
Indiana Univ. Math. J., {\bf 35} (1986), 245-267.

\bibitem[GL87]{gl87}
N. Garofalo and F.H. Lin,  \emph{Unique continuation for elliptic
operators: a geometric-variational approach}, Comm. Pure Appl.
Math., {\bf 40}, 347-366, 1987.

\bibitem[Gi83]{gi83}
M. Giaquinta, \emph{Multiple Integrals in The Calculus of Variations and Nonlinear Elliptic Systems}, Princeton University Press, Princeton, New Jersey, 1983.

\bibitem[GT83]{gt83}
D. Gilbarg and N. Trudinger, \emph{Elliptic Partial Differential Equations of Second Order}, 2nd Edition, Berlin Heidelberg New York 1983.


\bibitem[Jo13]{jo13}
J. Jost, \emph{Partial Differential Equations}, 3rd ed., Graduate Texts in Mathematics, Vol 214, Springer, 2013.


\bibitem[Ke07]{ke07}
C. Kenig, \emph{Some recent applications of unique continuation}, Contemp. Math., {\bf 439} (2007), 25-56.

\bibitem[KL88]{kl88}
V. A. Kondratiev and E. M. Landis,  \emph{Qualitative properties of the solutions of a second- order nonlinear equation}, Encyclopedia of Math. Sci. 32 (Partial Differential equations III), Springer-Verlag, Berlin (1988).

\bibitem[Ku98]{ku98}
I. Kukavica, \emph{Quantitative uniqueness for second order elliptic
operators}, Duke Math. J., {\bf 91} (1998), 225-240.

\bibitem[Me92]{m92}
V. Z. Meshkov, \emph{On the possible rate of decay at infinity of solutions of second order partial differential equations}, Math. USSR Sbornik, \textbf{72} (1992), 343-360.

\bibitem[LW13]{lw13}
C. L. Lin and J. N. Wang, \emph{Quantitative uniqueness estimates for the general second order elliptic equations},  J. Funct. Anal., \textbf{266} (2014), 5108-5125.

\bibitem[Sc98]{sc98}
F. Schulz, \emph{On the unique continuation property of elliptic divergence form equations in the plane}, Math. Z., \textbf{228} (1998), 201-206.

\bibitem[Ve62]{ve62}
I.N. Vekua, \emph{Generalized Analytic Functions}, Pergamon Press, London, 1962.

\end{thebibliography}
\end{document}